# How to Obtain Desirable Transfer Functions in MIMO Systems Using 2-Degrees of Freedom Controllers Under Internal Stability

Panos J. Antsaklis and Eloy Garcia

*Abstract*—This paper summarizes important results useful for controller design of multi-input multi-output systems. The main goal is to characterize the stabilizing controllers that provide a desired response.

## I. INTRODUCTION

THE aim of this brief paper is to use results from [1] to show how in Multi-Input Multi-Output (MIMO) systems, using 2-degrees of freedom compensation and existing techniques, the attainable response under internal stability can be completely characterized. These include the model matching problem, the inverse problem, the diagonal and block diagonal decoupling problems, and the static decoupling problem among others.

The requirement for stability restricts the design when the given system is non-minimum phase. Specifically, the compensated system transfer function matrix must have among its zeros all the unstable zeros of the given plant, together with their associated directions. In addition, the requirement for causality does not allow the compensated system transfer function to have relative degree lower than the given plant's.

This paper also shows constructively how to design feedback and feed-forward 2-degrees of freedom controllers that attain the desirable responses under internal stability. In fact all such controllers will be parametrically characterized. It is then up to the designer to decide how much of the needed control should be implemented via feedback compensation and how much via feedforward compensation. This decision depends on the specifications on feedback properties such as sensitivity to plant parameter variations, disturbance reduction, noise attenuation etc.

In many cases, the control architecture may not be a full 2-degrees of freedom architecture but may be restricted. For example it could be the very common in practice state feedback architecture, state feedback with observer, unity output feedback architecture etc. The effect of such restrictions on the attainable responses is examined and it is shown in each case how to derive the appropriate controllers.

Special attention is being paid to the special case of the transfer function of the plant $P(s)$ being square and non-singular and it is shown how expressions are simplified; note that this case includes the SISO case.

Examples are used throughout.

P. J. Antsaklis is with the Department of Electrical Engineering, University of Notre Dame, Notre Dame, IN 46556 USA.
E. Garcia is with the Control Science Center of Excellence, Air Force Research Laboratory, Wright-Patterson AFB, 45433.

In the following section, Section II, two basic results are first presented. The first characterizes all responses attainable using a 2-degrees of freedom controller under internal stability. The second result parametrically characterizes all stabilizing 2-degrees of freedom controllers, which attain the desired response, showing explicitly the interplay of feedforward and feedback control actions. Architectures that realize and implement the control policies are also shown.

In Section III several design problems are discussed including the decoupling and inverse problems.

The effect of several more restricted common control architectures on the system response is discussed in Section IV.

In Section V we study the special case of the plant transfer function matrix being square and nonsingular which includes the SISO case and derive results for the case when the denominator of the compensated transfer function $T$ is given.

## II. DESIGN OF CONTROL SYSTEMS FOR DESIRED RESPONSE

Consider the MIMO plant with $p \times m$ transfer function matrix $P = ND^{-1}$ ($y = Pu$), where $N$, $D$ are right coprime polynomial matrices. (The results in this paper are based on [1], pp. 622-644).
Suppose we want to control $P$ so that:

$$\begin{bmatrix} y \\ u \end{bmatrix} = \begin{bmatrix} T \\ M \end{bmatrix} r$$

where $T$, $M$ are proper and stable desired transfer function matrices and $r$ is an external input.

Basic Result I. It has been shown in [1] (Th. 4.23 – p. 627; also Th. 4.24 for proper stable factorizations) that the above $T$ and $M$ can be realized via a general 2-degrees of freedom controller with internal stability <u>if only if</u> there exist stable rational matrix $X$ so that:

$$\begin{bmatrix} T \\ M \end{bmatrix} = \begin{bmatrix} N \\ D \end{bmatrix} X.$$

In [1], Fig. 7.12 – p. 627, shows a way to realize such $T$ and $M$. See Fig. 1 below where $C_y$ is any feedback stabilizing controller for $P$.

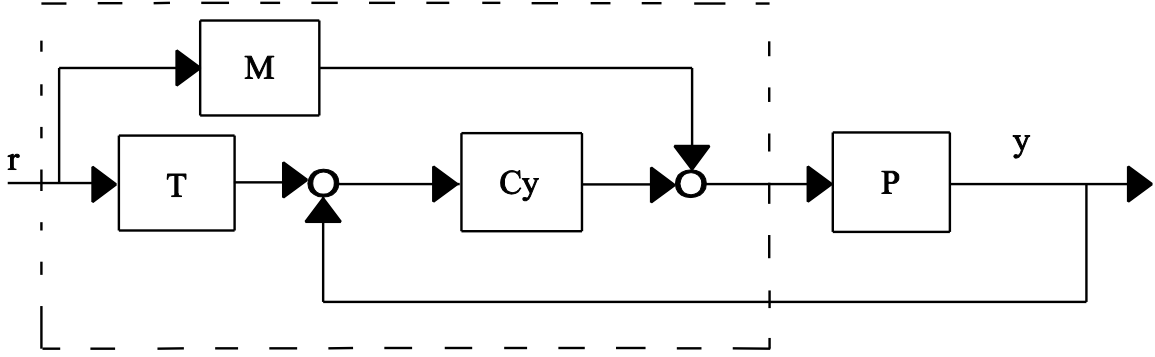

Fig. 1.

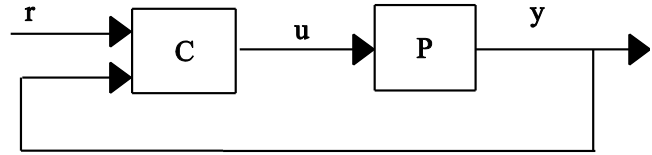

Fig. 2.

Basic Result II. The general 2-degree of freedom configuration is shown in Fig. 2, where $u = C \begin{bmatrix} y \\ r \end{bmatrix}$ = $[C_y, \; C_r] \begin{bmatrix} y \\ r \end{bmatrix}$ (see [1] p. 623 – Fig. 7.10 and Th. 4.21), with $C_y$ any feedback stabilizing controller ($u = C_y y$ internally stabilizes $y = Pu$) and $C_r$ such that $X = D^{-1}(I - C_y P)^{-1} C_r$ stable or $C_r = (I - C_y P) DX = (I - C_y P) M$.

Here $y = Tr = NXr$ and $u = Mr = DXr$.

*Feedback Controller* $C_y$. $C_y$ can be any stabilizing controller. To see the flexibility we characterize all stabilizing controllers ([1] Th. 4.22 – p. 624):

$$C = [C_y, \; C_r] = (I + QP)^{-1}[Q, \; M]$$
$$= [(I + LN)D^{-1}]^{-1}[L, \; X]$$

where $Q = DL$, $M = DX$ are proper with $L$, $X$, $D^{-1}(I + QP) = (I + LN)D^{-1}$ stable (also $(I + QP)^{-1}$ exists and is proper).

$$= [X, \; -KN]^{-1}[-(X_2 + KD), \; X]$$

where $K$, $X$ stable; $K$ is the Youla parameter. Other expressions may be found in [1] Th. 4.22 – p.624.

### III. RESPONSE.

In our case $X$ is given (or part of $X$ is specified). Note that $T$, $M$ cannot be chosen independently for a solution to (1) to exist, but they must be in the range of $\begin{bmatrix} N \\ D \end{bmatrix}$, or equivalently,

$\text{rank} \begin{bmatrix} N & T \\ D & M \end{bmatrix} = \text{rank} \begin{bmatrix} N \\ D \end{bmatrix} = m$, full column rank. They also must be stable, which is to be expected.

So, the only significant restrictions imposed (other than properness) are due to the unstable zeros of $P$, (in $N$), if any. Since $T = NX$ and $T$, $X$ stable, $T$ must contain <u>all</u> unstable zeros of the plant $P$.

The example 4.12 in [1], on p. 628 where $P = \dfrac{(s-1)(s+2)}{(s-2)^2}$ illustrates all these points.

It is also shown in the next section what additional restrictions are imposed if instead of a general 2-degrees of freedom configuration the more restricted unity feedback is used.

*Example*. We consider $P = \dfrac{(s-1)(s+2)}{(s-2)^2}$ and wish to characterize all proper and stable transfer functions $T(s)$ that can be realized by means of some control configuration with internal stability. Let $P = \dfrac{(s-1)}{(s+2)} \left[ \dfrac{(s-2)^2}{(s+2)^2} \right]^{-1} = N'D'^{-1}$ be an rc MFD in $RH_\infty$. Then in view of [1], Theorem 4.24 on p. 627, all such $T$ must satisfy

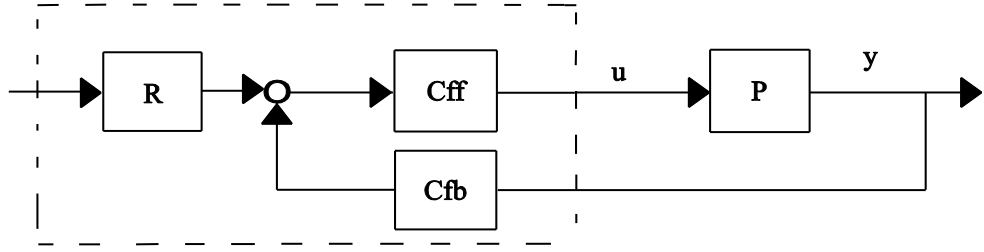

Fig. 3

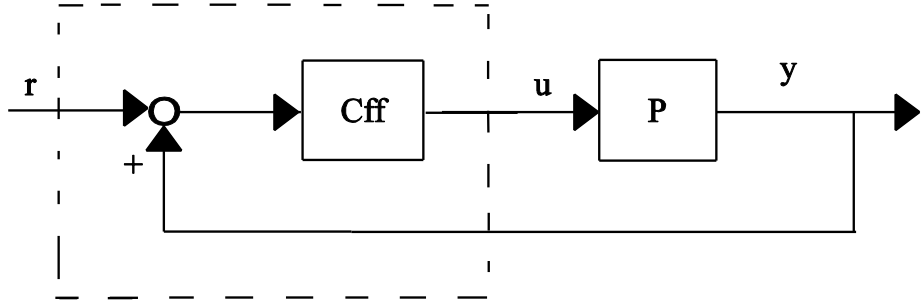

Fig. 4

$$N'^{-1}T = \frac{(s+2)}{(s-1)}T = X' \in RH_\infty$$

Therefore, any proper T with a zero at +1 can be realized via a 2-degrees of freedom feedback controller with internal stability.

Diagonal Decoupling. Here the desired transfer function $T(s)$ is not completely specified but is required to be diagonal, proper and stable. That is

$$N'X' = diag(n_i/d_i)$$

where $N'$, $X'$ are proper and stable. If $P^{-1}$ exists (see also the last section of this paper) then

$$X' = N'^{-1} diag(n_i/d_i).$$

If $P(s)$ has only stable zeros then no additional restrictions are imposed on $T(s)$.

When diagonal decoupling is accomplished using more restricted control configurations (e.g. unity feedback) additional restrictions are imposed on $X'$ due to the control architecture (see next section).

When linear state feedback $u = Fx + Gr$ is used for diagonal decoupling $T = N'D_f'^{-1}G$, $M = D'D_f'^{-1}G$, and $X' = D_f'^{-1}G$.

So, when $P^{-1}$ exists

$$X' = D_f'^{-1}G = N'^{-1}diag(n_i/d_i)$$

for a solution to exist (see Exercise 4.17, 4.20, Chapter 4 of [1]).

Inverse problem. In this case $T(s)=I$ and so

$$T = I = N'X'$$

which gives the conditions on the stable and proper parameter $X'$ (see also Exercise 4.17, 4.18 in Chapter 4 of [1]).

Static decoupling. Here $T(s)$ is square, proper and stable and $T(s)=\Lambda$, a real nonsingular diagonal matrix, so that a step change in the first input will only affect the first output at steady state. From

$$T = N'X'$$

in view of $T(s)=\Lambda$, a necessary and sufficient condition is that $P(s)$ does not have a zero at the origin (so $\det N'(0) \neq 0$). Then $X'(s)$ must satisfy $X'(s) = \left(N'^{-1}(s)\right)\Lambda$. This is also sufficient and static decoupling in this case can be achieved with just pre-compensation by a real gain matrix $C_r$

$$C_r = P^{-1}(s)\Lambda.$$

This can be seen from the expression in Section II for

$$C_r = (I - C_y P)D'X' = D'X' \quad (C_y = 0)$$

from which

$$C_r(s) = D'(s)X'(s) = P^{-1}(s)\Lambda.$$

## IV. RESTRICTED CONTROLLER CONFIGURATION.

If the configuration in Fig. 3 ([1], Fig. 7.13 – p. 629) is chosen, then

$$u = [C_y, \ C_r]\begin{bmatrix}y\\r\end{bmatrix} = [C_{ff}C_{fb}, \ C_{ff}R]\begin{bmatrix}y\\r\end{bmatrix} = D_c'^{-1}[N_y', \ N_r']$$

a lc factorization is $RH_\infty$ (proper and stable).

Then, if $C = [C_y, \ C_r]$ is given one could choose

$$R = N_r', \ C_{ff} = D_c'^{-1}, \ C_{fb} = N_y'.$$

Note that $R$, $C_{fb}$ are stable; $C_{ff}^{-1}$ exists and it is stable. In terms of the proper and stable parameters $X'$, $K'$ they can be related as

$$R = X', \ C_{ff} = (X', \ -K'N')^{-1}, \ C_{fb} = (X'_2 + K'N').$$

If the configuration in Fig. 4 (unity feedback) is chosen, ([1], Fig. 7.15 – p. 631), then

$$u = [C_y, \ C_r]\begin{bmatrix}y\\r\end{bmatrix} = [C_{ff}, \ C_{ff}]\begin{bmatrix}y\\r\end{bmatrix}$$

i.e. $C_r = C_y$ and in view of (4.158) (or 4.161)

$$L = X \ (L' = X') \text{ or } (Q = M).$$

Here (from 4.158)

$$C_{ff} = (I + MP)^{-1}M = M(I + PM)^{-1}$$
$$= M(I + T)^{-1} = [(I + XN)D^{-1}]^{-1}X$$

where $D^{-1}(I + MP) = (I + XN)D^{-1}$ stable which imposes a restriction on the allowed $X$ (or $X$').

Consider again the example above, if a single degree of freedom controller must be used, the class of realizable $T(s)$ under internal stability is restricted. In particular, if the unity feedback configuration $\{I; C_{ff}, I\}$ in Fig. 4 is used, then all proper and stable $T$ that are realizable under internal stability are again given by

$$T = N'X' = \frac{(s-1)}{(s+2)}X'$$

where $X' = L' \in RH_\infty$ and in addition

$$(I + X'N')D'^{-1} = \left[1 + X'\frac{(s-1)}{(s+2)}\right]\frac{(s+2)^2}{(s-2)^2} \in RH_\infty$$

i.e., $X' = n_x / d_x$ is proper and stable and should also satisfy

$$(s+2)d_x + (s-1)n_x = (s-2)^2 p(s)$$

for some polynomial $p(s)$.

This illustrates the restrictions imposed by the unity feedback controller, as opposed to a 2-degrees of freedom controller. Notice that these additional restrictions are imposed because the given plant has unstable poles.

The unity feedback case is studied in [1] Problem 7.23 p. 642-643.

## V. $P^{-1}$ EXISTS.

*Example.* Next it is shown how the expressions are simplified when $P^{-1}$ exists. This is from [1] Problem 7.23 – part (c) p. 643.

In the unity feedback configuration when $P^{-1}$ exists the closed loop system is internally stable if and only if:

$$[N^{-1}(I-T)P, \ N^{-1}T]$$

is stable and if $T$ is chosen so that

$$T = NX, \ X \text{ stable},$$

then this is the same as $(I + XN)D^{-1}$ stable (in Fig. 7.19 the sign is different).

In the Single-Input Single-Output (SISO) case, from (d) p. 643 when $P$ and $T$ are scalar, the conditions become

$$(1-T)d^{-1} = Sd^{-1} \text{ and } Tn^{-1} \text{ stable.}$$

The controller $C_{ff}$ is still given by the formula

$$C_{ff} = [(I + XN)D^{-1}]^{-1}X.$$

**Applications of results.**

Suppose then that the desired stable transfer function matrix $T$ needs to have only some desired denominator. In $T = N_T D_T^{-1}$, $D_T$ is given (note that in [1] p. 633, problems where $T$ is not specified completely, like in decoupling, are briefly discussed).

From $T = NX$, $X$ stable implies that we can select $X = D_T^{-1}$, assuming that $N_T = N$; other options exist where for example $X = N_x D_T^{-1}$ and $N_T = NN_X$.

From the previous discussion, for stability (in unity feedback) we need

$$X, \ (I + XN)D^{-1} \text{ stable, or}$$
$$D_T^{-1}, \ D_T^{-1}(D_T + N)D^{-1} \text{ stable}$$

and since $D_T^{-1}$ is assumed stable, the only condition is

$$(D_T + N)D^{-1} \text{ stable}$$

and the controller is given by

$$C_{ff} = [(I + D_T^{-1}N)D^{-1}]^{-1}D_T^{-1}$$
$$= (D_T^{-1}(D_T + N)D^{-1})^{-1}D_T^{-1}$$
$$= D(D_T + N)^{-1}D_T D_T^{-1}$$
$$= D(D_T + N)^{-1}.$$

Consider now Fig. 4 and verify response

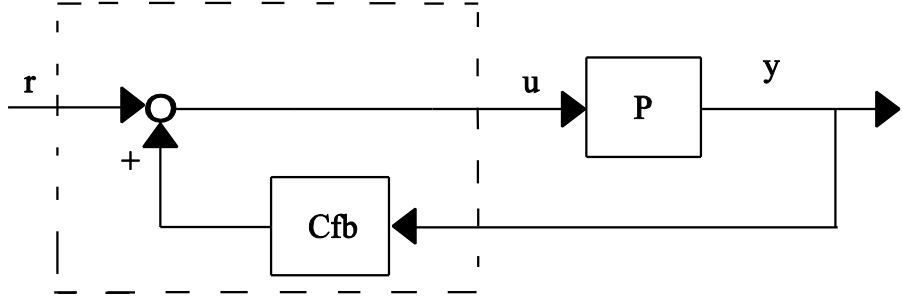

Fig. 5.

$$y = PC_{ff}(r+y) \rightarrow y = (I - PC_{ff})^{-1} PC_{ff} r.$$

Here $P = ND^{-1}$, $C_{ff} = D(D_T + N)^{-1}$

$$\begin{aligned} y &= PC_{ff}(I - PC_{ff})^{-1} r \\ &= ND^{-1}D(D_T + N)^{-1}\left[I - ND^{-1}D(D_T + N)^{-1}\right]^{-1} r \\ &= N(D_T + N)^{-1}(D_T + N)[D_T + N - N]^{-1} r \\ &= ND_T^{-1} r \end{aligned}$$

as expected.

Note that here

$$T^{-1} + P^{-1} = D_T N^{-1} + DN^{-1} = (D_T + D)N^{-1}$$
$$= \left[(D_T + D)D^{-1}\right]\left[DN^{-1}\right] = C_{ff}^{-1} P^{-1}.$$

If the configuration in Fig. 5 is chosen (see also [1], Fig. 7.18 – p. 632), then

$$u = [C_y, \ C_r]\begin{bmatrix} y \\ r \end{bmatrix} = [C_{fb}, \ I]\begin{bmatrix} y \\ r \end{bmatrix}$$

and then (4.158) – p.624

$$C_{fb} = \left[(I + LN)D^{-1}\right]^{-1} L$$

where $L$ is restricted to satisfy

$$X = (I + LN)D^{-1} \qquad (XD = I + LN).$$

In the special case when $P^{-1}$ exists ($T^{-1}$ exists)

$$C_{fb} = X^{-1}(XD - I)N^{-1} = X^{-1} L.$$

So the stability condition is just $X$ must be stable. Here,

$$\begin{aligned} y &= (I - PC_{fb})^{-1} Pr = P(I - C_{fb} P)^{-1} r \\ &= (ND^{-1})\left[I - [(I + LN)D^{-1}]^{-1} LND^{-1}\right]^{-1} r \\ &= (ND^{-1})\left[((I + LN)D^{-1})^{-1}[(I + LN) - LN]D^{-1}\right]^{-1} r \\ &= (ND^{-1})\underbrace{(I + LN)D^{-1}}_{X} r = NXr \end{aligned}$$

as desired.

Assuming $P^{-1}$ and $T^{-1}$ exist, when $X = D_T^{-1}$ stable form previous page, for internal stability

$$\begin{aligned} C_{fb} &= (D_T^{-1})^{-1}\left[D_T^{-1} D - I\right] N^{-1} \\ &= D_T D_T^{-1}[D - D_T] N^{-1} \\ \Rightarrow C_{fb} &= (D - D_T) N^{-1}. \end{aligned}$$

Clearly, the closed loop transfer function with $C_{fb}$ as above

$$\begin{aligned} T &= P(I - C_{fb} P)^{-1} = (ND^{-1})\left[I - (D - D_T)N^{-1} ND^{-1}\right]^{-1} \\ &= (ND^{-1})D\left[D - D + D_T\right]^{-1} \\ &= ND_T^{-1} \end{aligned}$$

as expected.

So given $T = NX = ND_T^{-1}$ stable ($D_T$ desired) select

$$C_{fb} = (D - D_T) N^{-1}$$

and only worry about properness. Note that here

$$T^{-1} - P^{-1} = D_T N^{-1} - DN^{-1} = (D_T - D)N^{-1} = -C_{fb}.$$

## VI. Concluding remarks

Several results from [1] on 2-degree of freedom controllers were used to solve a set of problems where the denominator of the desired transfer function $T$ is given. The control design problems where $T$ or parts of it are specified can be addressed in similar way. Different control configurations can also be studied; for example, in the linear state feedback case $u = Fx + Gr$, $X = (D_f)^{-1} G$.

## References

[1] P. J. Antsaklis and A. N. Michel, *Linear Systems*, Springer/Birkhauser, 2006.